\documentclass[pdflatex]{sn-jnl}
\usepackage{enumitem}
\usepackage{graphicx}%
\usepackage{anyfontsize}
\usepackage{lmodern}
\usepackage{multirow}%
\usepackage{amsmath,amssymb,amsfonts}%
\usepackage{amsthm}%
\usepackage{mathrsfs}%
\usepackage[title]{appendix}%
\usepackage{xcolor}%
\usepackage{textcomp}%
\usepackage{manyfoot}%
\usepackage{booktabs}%
\usepackage{algorithm}%
\usepackage{algorithmicx}%
\usepackage{algpseudocode}%
\usepackage{listings}%
\usepackage[utf8]{inputenc}
\usepackage{derivative}
\usepackage{epstopdf}
\usepackage[caption=false]{subfig}
\usepackage{mathtools}
\usepackage{nicefrac}
\usepackage[numbers,sort&compress]{natbib}% Citation support using natbib.sty
\bibpunct[, ]{[}{]}{,}{n}{,}{,}% Citation support using natbib.sty
\makeatletter% @ becomes a letter
\def\NAT@def@citea{\def\@citea{\NAT@separator}}% Suppress spaces between citations using natbib.sty
\makeatother% @ becomes a symbol again
\usepackage{multirow}

\usepackage{microtype}
\theoremstyle{plain}% Theorem-like structures provided by amsthm.sty
\newtheorem{theorem}{Theorem}[section]

\newtheorem{assumption}[theorem]{Assumption}
\usepackage{algorithm}
\usepackage{algcompatible}
\usepackage{csvsimple}
\usepackage{graphicx}
\usepackage{algpseudocode}
\theoremstyle{definition}

\theoremstyle{remark}

\DeclareMathOperator{\tridiag}{tridiag}
\usepackage{subfig}
\usepackage{pgfplots,tikz}
\usetikzlibrary{shapes.geometric}
\pgfplotsset{compat=1.18}

\usepackage{booktabs}
\usepackage{multirow}
\newcommand{\mc}{\mathrm{i}}
\begin{document}
%%%%%%%%%%%%%%%%%%%%%%%%%%%%%%%%%%%
\title{ Scaling-and-squaring method for computing the inverses of matrix $\varphi$-functions}

\author[1]{\fnm{Lidia} \sur{Aceto}}\email{lidia.aceto@uniupo.it}

\author*[2]{\fnm{Luca} \sur{Gemignani}}\email{luca.gemignani@unipi.it}
%\equalcont{These authors contributed equally to this work.}

\affil[1]{\orgdiv{Dipartimento di Scienze e Innovazione Tecnologica}, \orgname{Universit\`a del Piemonte Orientale}, \orgaddress{\street{viale Teresa Michel, 11}, \city{Alessandria}, \country{Italy}}}

\affil*[2]{\orgdiv{Dipartimento di Informatica}, \orgname{Universit\`a di Pisa}, \orgaddress{\street{Largo Bruno Pontecorvo, 3}, \city{Pisa}, \country{Italy}}}

%\begin{document}

%%%%%%%%%%%%%%%%%%%%%%%%%%%%%%%%%%%%%%%%%%%%%%%%%%%%%%%%%
\abstract{This paper aims to develop efficient numerical methods for computing the inverse of matrix $\varphi$-functions, $\psi_\ell(A) := (\varphi_\ell(A))^{-1}$,  for $\ell =1,2,\ldots,$ when $A$ is a large and sparse matrix with eigenvalues in the  open left half-plane. While $\varphi$-functions play a crucial role in the analysis and implementation of exponential integrators, their inverses arise in solving certain direct and inverse differential problems with non-local boundary conditions. We propose an adaptation of the standard scaling-and-squaring technique for computing $\psi_\ell(A)$, based on the Newton-Schulz iteration for matrix inversion. The convergence of this method is analyzed  both theoretically and numerically. In addition, we derive and analyze Padé approximants for approximating $\psi_1(A/2^s)$, where $s$ is a suitably chosen integer, necessary at the root of the squaring process. Numerical experiments demonstrate the effectiveness of the proposed approach.}


\maketitle
%%%%%%%%%%%%%%%%%%%%%%%%%%%%%%%%%%%%%%%%%%%%%%%%%%%%%%%%%
\section{Introduction}

In many applications involving complex systems with widely varying time scales, exponential integrators play a central role in computing dynamics. They require the efficient and accurate evaluation of matrix functions, which are closely related to the matrix exponential and are referred to as $\varphi$-functions in the recent literature \cite{HLS,EXPINT}.  The $\varphi$-functions are defined as follows
\begin{eqnarray} \label{integrale}
 \varphi_0(z)=e^z,   \qquad \varphi_\ell(z) = \frac{1}{(\ell-1)!} \int_0^1 e^{(1-\xi)z} \xi^{\ell-1} \operatorname{d\xi},  \ \ell\geq 1, \qquad  z \in \mathbb{C},
\end{eqnarray}
where, for integers $\ell \ge 1,$ the second expression provides the integral representation of $\varphi_\ell(z).$
These  $\varphi$-functions  are entire functions with  the Taylor series expansion
\begin{eqnarray} \label{Taylor}
\varphi_{\ell}(z)=\sum_{k=0}^{+\infty} \frac{z^k}{(\ell+k)!}= \left( e^z- \displaystyle\sum_{j=0}^{\ell-1} \displaystyle\frac{z^j}{j!}\right)z^{-\ell}, \quad  \ell\geq 0.
\end{eqnarray}
This expansion can be extended to a matrix argument by setting 
\[
\varphi_{\ell}(A)=\sum_{k=0}^{+\infty} \frac{A^k}{(\ell+k)!}, \quad  \ell\geq 0, \quad A\in \mathbb R^{n\times n}.
\]
The reciprocal of the $\varphi$-functions, referred to as the $\psi$-functions, is given by 
\[
\psi_\ell(z):=\frac{1}{\varphi_\ell(z)}= \left( e^z- \displaystyle\sum_{j=0}^{\ell-1} \displaystyle\frac{z^j}{j!}\right)^{-1} z^\ell, \quad  \ell\geq 0.
\]
For small values of $\ell,$  these functions take the following form 
\begin{align*}
\psi_0(z)=e^{-z},  \quad \psi_1(z)=\displaystyle\frac{z}{e^z-1}, \quad \psi_2(z)=\displaystyle\frac{z^2}{e^z-1 -z }.   
\end{align*}
The $\psi$-functions are meromorphic functions with  interesting applications in the solution of  direct and inverse  boundary value problems  for  abstract first order differential systems with non-local boundary conditions \cite{Boito,DB,aceto2024,PKT,KT}. These problems also  arise from parabolic equations  that are discretized  by using a semi-discretization in space (method of lines) \cite{PRIT,PKK}.

In this paper we focus on developing  efficient numerical methods  for computing 
$\psi_\ell(A)$.   Specifically, we propose a scaling-and-squaring  algorithm tailored for matrices $A$ with eigenvalues that have negative real parts. The  Newton-Schulz iteration for  matrix inversion  given in \cite{PR,Pan} can be regarded as  the 
work-horse of our scheme.  This iteration is quadratically convergent and only requires matrix-by-matrix products and, hence, it can be implemented with great efficiency on systolic arrays and parallel computers. A careful analysis of this iteration provides convergence results based on the properties of the scalar $\varphi_\ell(z)$ and $\psi_\ell(z)$ functions involved.  
Additionally, we exploit rational approximations of these functions to determine appropriate initial guesses for the iteration.

A  scaling-and-squaring method for computing $\varphi_\ell(A)$ 
was proposed in \cite{Skaflestad},  which efficiently updates 
$\varphi_\ell(A) $ to $\varphi_\ell(2 A)$ using recurrence relations. However, a direct analogous approach cannot be applied to $\psi$-functions.  Instead, we show  that $\psi_\ell(2 A)$ can be computed iteratively, starting with $\varphi_\ell(2A)$  and $\psi_\ell(A),$  by using the  Newton-Schulz iteration  for matrix inversion. The initial guess for this iteration is $\psi_\ell(A),$ and the matrix inversion is applied to $\varphi_\ell(2 A).$  
The convergence condition for this iterative method is
\[
\left|1- \displaystyle\frac{\varphi_\ell(2 \lambda)}{\varphi_\ell(\lambda)}\right|< 1,
\]
which holds when all eigenvalues $\lambda$ of the matrix $A$ lie in the open left half-plane.

Our proposed scaling and modified squaring method for
evaluating the  matrix $\psi$-functions combines the squaring scheme based on the  Newton-Schulz iteration with  an efficient approximation of $\psi_\ell(A/2^{s})$, for an integer $s$ large enough so that  $A/2^{s}$ has small eigenvalues. Assuming these  eigenvalues lie within the strip $\mathbb  R \times \mc [-1, 1]$, we leverage results from \cite{gem2023} to conclude that the Newton-Schulz iteration, starting from 
$\psi_1(A/2^{s}),$ can be employed efficiently to compute 
$\psi_\ell(A/2^{s}).$

Functional approximation techniques are used to evaluate the starting point $\psi_1(A/2^{s})$. While rational approximations  of $\psi_1(z)$ have been derived using Fourier theory in \cite{Boito,Boito1,Boito3}, this paper presents an alternative approach based on Padé approximation. This method is particularly effective for deriving error estimates in a small region around the origin in the complex plane. We compute the diagonal Padé approximant for  $\psi_1(z)$ and provide a detailed analysis of the corresponding approximation error.  

By combining the scaling-and-squaring method with Padé approximation, we propose a feasible and efficient approach for computing $\psi_\ell(A)$. Several numerical experiments are presented to validate the effectiveness of this approach.

The paper is structured as follows. In Section~\ref{sec:invfipade},  we describe 
the construction of diagonal  Padé approximants for $\psi_\ell(z)$ with emphasis on the case $\ell=1$.  Section~\ref{sec:scalesquar}  presents  our scaling-and-squaring method for computing $\psi_\ell(A)$.  In Section~\ref{numres}, we provide numerical results that validate the effectiveness of our approach. Finally, Section~\ref{sec:fine} offers conclusions and outlines potential directions for future work.

%%%%%%%%%%%%%%%%%%%%%%%%%%%%%%%%%%%%%%%%%%%%%%%%%%%%%%%%%%%%%%%%%%%%%%%%%%%%
\section{Pad\'e  approximation of  \texorpdfstring{$\psi_\ell(z)$}{Lg}} \label{sec:invfipade}

In addition to the integral representation in \eqref{integrale}, 
 the $\varphi$-functions can also be expressed as special cases of confluent hypergeometric functions. It is known that the confluent hypergeometric function of the first kind is defined by
\begin{equation} \label{defconfhyp}
{}_{1}F_{1}(a;b;z) =\sum _{{k=0}}^{+ \infty}{\frac {a^{{(k)}}z^{k}}{b^{{(k)}}k!}}, \quad z \in \mathbb{C}, a \in \mathbb{C}, b \in \mathbb{C} \setminus \mathbb{Z}_0^-,
\end{equation}
where
\begin{equation} \label{pocham}
a^{(0)}=1, \, a^{{(k)}}=a(a+1)(a+2)\cdots (a+k-1)\,,
\end{equation}
denotes the rising factorial.
From this definition and the properties of the rising factorial \eqref{pocham}, we can derive the following relation (see \eqref{Taylor})
\begin{equation} \label{conflphi}
 \frac{1}{\ell!} {}_{1}F_{1}(1;\ell+1;z) = \frac{1}{\ell!}\sum_{k=0}^{+\infty} \frac{z^k}{(\ell+1)^{(k)}} =
 \sum_{k=0}^{+\infty} \frac{z^k}{(\ell+k)!} =\varphi_\ell(z).
\end{equation}
Given this fact, the general form of the $[d/d]$-Padé approximant $\widehat R_{d,\ell}(z)$ for $\varphi_\ell(z)$
 can be determined using either the results from Luke's work \cite{Luke}, which address the diagonal Padé approximant of $ {}_{1}F_{1}\left(1;b ;-z \right),$ or the results provided in \cite[Lemma 2]{Skaflestad} for $\varphi_\ell(z).$ In both cases, the approximant takes the form
\begin{equation} \label{padephi}
    \widehat R_{d,\ell}(z) = \frac{\widehat N_{d,\ell}(z)}{\widehat D_{d,\ell}(z)},
\end{equation}
where the polynomials $\widehat N_{d,\ell}(z)$ and $\widehat D_{d,\ell}(z)$ are
\begin{eqnarray} 
   \widehat N_{d,\ell}(z) &=& \frac{d!}{(2d+\ell)!} \sum_{i=0}^d \left[ \sum_{j=0}^i \frac{(2d+\ell-j)! (-1)^j}{j!(d-j)! (\ell+i-j)!}\right] z^i, \label{main1}\\
   \widehat D_{d,\ell}(z) &=& \frac{d!}{(2d+\ell)!} \sum_{i=0}^d \frac{(2d+\ell-i)!}{i!(d-i)! } (-z)^i. \label{main2}
\end{eqnarray}

The difference lies in the accuracy of the error estimates: the one derived from the diagonal Padé approximant of $\varphi_\ell(z),$ as presented in \cite[Lemma 2]{Skaflestad}, is sharper than the estimate provided in \cite{Luke}.
Therefore, we report here the more accurate  asymptotic error estimates given in \cite{Skaflestad}, namely 
\begin{equation*} \label{stimaFI0}
    \varphi_\ell (z) - \widehat R_{d,\ell}(z) = \frac{(-1)^d d! (d+\ell)!}{(2d+\ell)! (2d+\ell+1)!} z^{2d+1} + \mathcal{O}(z^{2d+2}).
\end{equation*}
A precise non-asymptotic error estimate is useful for obtaining an upper bound on the approximation error. This estimate can be derived by applying the Cauchy product to the power series expansions of $\varphi_\ell(z)$ (as given in \eqref{conflphi}) and the denominator $\widehat D_{d,\ell}(z)$. Specifically, we find that 
\begin{eqnarray*} 
\varphi_\ell (z)\widehat D_{d,\ell}(z) 
     &=& \sum_{i=0}^{+\infty} \sum_{j=0}^{\min\{i,d\}} \frac{ d! (2d+\ell-j)! (-1)^j}{(2d+\ell)! j! (d-j)! (\ell+i-j)!} z^{i}\\
     &=&\widehat N_{d,\ell}(z)  +
    \sum_{i=d+1}^{+\infty} \sum_{j=0}^{d}\frac{ d! (2d+\ell-j)! (-1)^j }{ (2d+\ell)! j! (d-j)! (\ell+i-j)! } z^{i}.
\end{eqnarray*}
After some computation, using the identity
\[
\binom{n}{k} = (-1)^k \binom{k-n-1}{k}
\]
we obtain the relation
\begin{eqnarray*} 
\varphi_\ell (z)\widehat D_{d,\ell}(z) 
     &=&\widehat N_{d,\ell}(z)  +
     \sum_{i=d+1}^{+\infty} 
    \frac{(-1)^d d! (d+\ell)!}{(2d+\ell)!(\ell+i)!}
    \sum_{j=0}^{d} \binom{-d-\ell-1}{d-j} \binom{\ell+i}{j} z^{i}.
\end{eqnarray*}
Now, observe that the sum
$
\sum_{j=0}^{d} \binom{-d-\ell-1}{d-j} \binom{\ell+i}{j}
$
is the coefficient of $x^d$ in the expansion of  $(1+x)^{-d-\ell-1} (1+x)^{\ell+i} = (1+x)^{i-(d+1)}.$ Hence,  
\begin{itemize}
    \item for $(d+1)\leq i\leq 2d,$ 
    \[
    \sum_{j=0}^{d} \binom{-d-\ell-1}{d-j} \binom{\ell+i}{j}=0,
    \]
    \item for $i\geq (2d+1),$ 
    \[
    \sum_{j=0}^{d} \binom{-d-\ell-1}{d-j} \binom{\ell+i}{j}=\binom{i-(d+1)}{d}.
    \]
\end{itemize}
Thus, we arrive at the following relation
\begin{eqnarray} \label{stimaFI}
    \varphi_\ell (z)\widehat D_{d,\ell}(z) 
    - \widehat N_{d,\ell}(z) =\sum_{i=2d+1}^{+\infty}\frac{(-1)^d d! (d+\ell)!}{(2d+\ell)! (\ell+i)!} \binom{i-(d+1)}{d}z^{i},
    %=\frac{d!}{(2d+\ell)!}\sum_{i=2d+1}^{+\infty} \sum_{j=0}^{i}\frac{  (2d+\ell-j)! (-1)^j }{ j! (d-j)! (\ell+i-j)! } z^{i},
\end{eqnarray}
which can be used to derive upper bounds for the approximation error.\\

The diagonal Padé approximants for the $\psi$-functions, along with suitable error estimates, can be derived from \eqref{main1}--\eqref{stimaFI}  using the {\it reciprocal covariance} property \cite[p. 61]{Cuyt}, which we recall below.
\begin{theorem} \label{teo2.1}
    Let $r_{\mu,\nu}(z)= p_{\mu}(z)/q_{\nu}(z)$ %with $q_\nu(0)=1$ is 
    be the $[\mu/\nu]$-Padé approximant  for a function $f(z)$ with the formal power series  $f(z) = \sum_{j=0}^{+ \infty} c_jz^j,$ $c_j \in \mathbb{C},$   $ c_0 \neq 0.$ Then,  the $[\nu/\mu]$-Padé approximant for $1/f(z)$ is given by
\[
r_{\nu,\mu}(z)=\displaystyle\frac{q_{\nu}(z)/c_0}{p_{\mu}(z)/c_0}.
\]
\end{theorem}
In particular, from \eqref{conflphi} we known that in our case the coefficient $c_0=1/\ell!.$ Using this result and the expression from
 \eqref{padephi}, we can obtain the $[d/d]$-Padé approximant of $\psi_\ell(z)$:
\begin{equation} \label{main}
  \psi_\ell(z) \approx \frac{\widehat D_{d,\ell}(z) \ell!}{ \widehat N_{d,\ell}(z) \ell!} 
 = \frac{{\mathcal{N}}_{d,\ell} (z) }{{\mathcal{D}}_{d,\ell}(z)}.
\end{equation}

%%%%%%%%%%%%%%%%%%%%%%%%%%%%%%%%%%%%%%%%%%
\subsection{Exploring the case  \texorpdfstring{$\ell=1$}{Lg1}}

It is straightforward to derive the $[d/d]$-Padé approximant of $\psi_1(z)$  from \eqref{main} using the relationships in \eqref{main1} and \eqref{main2}. The resulting approximant is given by
\begin{equation} \label{padepsi1}
\mathcal{R}_{d, 1}(z) = \frac{{\mathcal{N}}_{d,1} (z)}{{\mathcal{D}}_{d,1}(z)}
\end{equation}
where the numerator and denominator are expressed as follows:
\begin{equation}\label{padepol}
\begin{split}
    {\mathcal{N}}_{d,1} (z) = & \; \frac{d!}{(2d+1)!} \sum_{i=0}^d \frac{(2d+1-i)!}{i!(d-i)! } (-z)^i,\\
    {\mathcal{D}}_{d,1} (z) = & \; \frac{d!}{(2d+1)!} \sum_{i=0}^d \left[ \sum_{j=0}^i \frac{(2d+1-j)! (-1)^j}{j!(d-j)! (1+i-j)!}\right] z^i. 
\end{split}
\end{equation}
We now analyze the error introduced when approximating 
$\psi_1(z)$ using its $[d/d]$-Padé approximant. Specifically, we 
derive an error estimate for this Padé approximation, starting from the relation given in \eqref{stimaFI}. To do so, we first recall the generating function for Bernoulli polynomials $B_j(\cdot)$, which is given by
\[
 \frac{z e^{zt}}{e^z-1} = \sum_{j=0}^{+\infty} B_j \left(t \right) \frac{z^j}{j!}, \qquad |z|<2 \pi,
\]
cf. \cite[Eq. 23.1.1]{Abramowitz}. 
The Bernoulli numbers $B_j$ are defined as $B_j = B_j(0),$ the values of the Bernoulli polynomials at $t=0.$ Using this generating function, we can readily derive the Taylor series expansion  for $\psi_1(z)$ as
\begin{equation} \label{seriePsi}
\psi_1(z) = \sum_{j=0}^{+\infty} \frac{B_j}{j!} z^j
\end{equation}
which converges for $|z|<2 \pi.$ 
Dividing both sides of \eqref{stimaFI} by $\varphi_\ell(z)$ and using \eqref{main}, we obtain
\begin{eqnarray*} 
   \mathcal{N}_{d,\ell}(z) 
    - \frac{1}{\varphi_\ell(z)} \mathcal{D}_{d,\ell}(z) = \frac{\ell!}{\varphi_\ell(z)}
    \sum_{i=2d+1}^{+\infty}\frac{(-1)^d d! (d+\ell)!}{(2d+\ell)! (\ell+i)!} \binom{i-(d+1)}{d}z^{i}.
\end{eqnarray*}
Setting $\ell=1$ and recalling that $1/\varphi_1(z)=\psi_1(z)$ 
\[
\psi_1(z){\mathcal{D}}_{d, 1}(z) -\mathcal{N}_{1,\ell}(z)  =  \psi_1(z) \sum_{i=2d+1}^{+\infty}\frac{(-1)^{d+1} d! (d+1)!}{(2d+1)! (1+i)!} \binom{i-(d+1)}{d}z^{i}.
\]
Using \eqref{seriePsi} and the Cauchy product of two series, we can rewrite the right-hand side of the previous formula as follows
\[
\begin{split}
\psi_1(z){\mathcal{D}}_{d, 1}(z) - \mathcal{N}_{d,1}(z) = 
\sum_{i=2d+1}^{+\infty} \sum_{m=0}^{i-2d-1}  \frac{B_m}{m!}   
  \frac{(-1)^{d+1} d!}{(2d+1)!}\frac{ (d+1)! }{ (1+i-m)!} \binom{i-m - (d+1)}{d} z^i.
\end{split}
\]
These considerations are summarized in the following result.
 \begin{theorem}\label{phipadetheo}
Let $d\geq 0$  be an integer, and let $z \in \mathbb C$ with  $|z|\le R <2 \pi.$ 
Define the $[d/d]$-Padé approximant for $\psi_1(z)$ as
\[
\mathcal{R}_{d, 1}(z) = \frac{{\mathcal{N}}_{d,1} (z)}{{\mathcal{D}}_{d,1}(z)}
\]
where the polynomials $\mathcal{N}_{d,1} (z)$ and $\mathcal{D}_{d,1} (z)$ are given by \eqref{padepol}.
Then, the following bound holds 
\begin{equation}\label{bb}
|\psi_1(z){\mathcal{D}}_{d, 1}(z)  -{\mathcal{N}}_{d,1} (z)| \leq  
\frac{ d! }{(2d+1)!} \sum_{i=2d +1}^{+\infty} \gamma_i  R^i =: s(d, R), \qquad |z|\leq R <2 \pi,
\end{equation}
where
\[
\gamma_i= \left|
\sum_{m=0}^{i-2d-1} \frac{B_m (d+1)!}{ m!(1+i-m)!} \binom{i-m-(d+1)}{d}
\right|.
\]
\end{theorem}

This theorem can be used to provide error estimates. For instance, for $d=6$ and $R=4$ numerical computations yield $s(6,4)\leq 9.8\cdot 10^{-7}$ and $\min\limits_{|z|\leq 4}|{\mathcal{D}}_{6, 1}(z)|\geq 0.53$ which together provide the bound 
\[
|\psi_1(z)  - \mathcal{R}_{6,1} (z) | \leq 1.9 \cdot 10^{-6}  \qquad \mbox{for } \, |z|\leq 4. 
\]
Similarly, for $d=9$ and $R=4$ we obtain 
$s(9,4)\leq 2.7\cdot 10^{-12}$ and $\min\limits_{|z|\leq 4}|{\mathcal{D}}_{9, 1}(z)|\geq 0.5$  which  give  the bound  
\[
|\psi_1(z)  - \mathcal{R}_{9,1} (z) | \leq 5.4 \cdot 10^{-12} \qquad \mbox{for } \, |z|\leq 4. 
\]

%%%%%%%%%%%%%%%%%%%%%%%%%%%%%%%%%%%%%%%%%%%%%%%%%%%%
%%%%%%%%%%%%%%%%%%%%%%%%%%%%%%%%%%%%%%%%%%%%%%%%%%%%
\section{The scaling-and-squaring method} \label{sec:scalesquar}

Building on the previous results, we propose a scheme for evaluating matrix $\psi$-functions. The three key components of this scheme are as follows:
\begin{itemize}[leftmargin=30pt] 
    \item an efficient method for evaluating  matrix $\varphi$-functions; 
    \item an efficient method for evaluating the matrix function $\psi_1(B),$ where $B$ is a matrix with eigenvalues clustered around the origin in the complex plane; 
    \item the Newton-Schulz iteration for matrix inversion.
\end{itemize}
The first key component is the scaling-and-squaring method, proposed in \cite{Skaflestad,EXPINT}, which is based on the recurrence relation 
\begin{equation} \label{squarfi}
 \varphi_\ell(2z) = \frac{1}{2^\ell} \bigg[ \varphi_0(z) \varphi_\ell(z) + \sum_{j=1}^\ell \frac{1}{(\ell-j)!} \varphi_j(z) \bigg], \quad \ell\geq 0,
\end{equation}
and is implemented in the algorithm {\tt phipade} of the package EXPINT. 

Assume that the input matrix $A$ is scaled such that $\| A/2^s \|_\infty \le \theta,$
for some positive constant $\theta$. The $[d/d]$-Padé approximations $\widehat R_{d,j}(z)\approx \varphi_j(z), j=0,1,\dots, \ell,$ are then computed with  $z=  A/2^{s}$ where 
\begin{equation}\label{scriterio}
 s=\operatorname{max}\bigg\{ \bigg\lceil
 \operatorname{log}_2\left( \frac{\|A\|_\infty}{\theta} \right)\bigg\rceil, 0 \bigg\}.
\end{equation}
Once the Padé approximations are computed,  the scaling is undone by performing $s$ iterations of the recurrence relation in \eqref{squarfi} to approximate $\varphi_\ell(A)$. 

As for the second tool, if  $\theta$ in \eqref{scriterio} is properly chosen, then $B=A/2^s$ will have all its eigenvalues within the convergence region specified in Theorem~\ref{phipadetheo}. Therefore, we can approximate 
 $\psi_1(A/2^s)$ using the Pad\'e approximants $\mathcal{R}_{d,1}(A/2^s)$.  According to  Theorem~\ref{teo2.1}, these  approximants  can be generated by simply swapping the numerator with the denominator of the approximations of $\varphi_1(A/2^s).$ This procedure is implemented in the algorithm {\tt psi1eval}, which takes  a matrix $B$ as input and returns an approximation of $\psi_1(B)$ as output, based on the method described in \cite {Skaflestad}.
  
The third  key ingredient is the Newton-Schulz iteration, which is used to invert a nonsingular  matrix $M\in \mathbb C^{n\times n}.$ The iteration is defined as
\begin{equation}\label{Newton}
X_0\in \mathbb C^{n\times n}, \qquad X_{k+1}=2 X_k  - X_k  M   X_k , \,\, k\geq 0.
\end{equation}
From the relation
\begin{equation}\label{quaderror}
I-X_{k+1} M= (I-X_k M)^2=(I-X_0 M)^{2^k},
\end{equation}
we can conclude that   the  Newton-Schulz iteration \eqref{Newton}   converges quadratically to $M^{-1},$ provided that all eigenvalues of the matrix $(I-X_0 M)$ have modulus less than $1.$ Furthermore,  from \eqref{quaderror} we deduce that 
\begin{equation}\label{stop_criterion}
\| X_{k+1}-X_k\|_{\infty} \leq \| X_{k}-X_{k-1}\|_{\infty}  \|(I-X_0 M)^{2^{k-1}}\|_{\infty}  \| X_{k-1}^{-1} X_k\|_{\infty }.
\end{equation}
This inequality indicates that the sequence  $\| X_{k}-X_{k-1}\|_{\infty}$ is expected to decrease monotonically  in the regime.   We exploit this property to develop an implementation of the Newton-Schulz iteration which does not depend on any input tolerance parameter. A basic implementation of the Newton-Schulz iteration is given in Algorithm~1.
%%%%%%%%%%%%%%%%%%%%%%%%%%%%%%%%%%%%%%%%%
\begin{algorithm}\caption{: Procedure {\textsc{NewtonSchulz}$(M, X_0)$}. Given the in\-vertible matrix $M$ and a suitable initial guess $X_0$  this procedure returns the approximation $X$ of $M^{-1}$ computed by the  Newton-Schulz iteration}
    \label{algorithm2}
    \begin{algorithmic}[1] 
        \State Initialization: $err_1\gets {\tt inf}$; $err_2\gets {\tt inf};$
            \WHILE {$err_1\geq err_2  \ | \ err_1>=0.1$} 
            \STATE $X \gets 2 X_0 - X_0 M X_0$; 
            \STATE $err_1\gets err_2$; $ err_2 \gets \| X-X_0 \|_\infty $; $X_0\gets X$; 
            \ENDWHILE
            \STATE \textbf{return} $X \approx M^{-1}$
    \end{algorithmic}
\end{algorithm}
%%%%%%%%%%%%%%%%%%%%%%%%%%%%%%%%%%%%%%%%%
It should be noted that each iteration of Procedure {\textsc{NewtonSchulz}} requires two matrix-by-matrix multiplications, making it particularly well-suited for applications involving structured matrices in high-performance computing environments.\\

These three components are now combined in the development of our proposed algorithm for computing matrix $\psi$-functions. Suppose the following condition holds:
\smallskip

\begin{assumption}\label{ass1} The  matrix $A$  has  all eigenvalues with negative real part. 
\end{assumption}
\smallskip
As will be shown in Subsection~\ref{convan}, under this assumption, we can apply the Newton-Schulz method to compute 
 $(\varphi_\ell(2 B))^{-1}=\psi_\ell(2B)$  starting from $\psi_\ell(B)$.  
 This process can be iterated at any squaring step, providing an effective way to approximate $\psi_\ell(A)$.  The resulting procedure, called {\textsc{UpdateMatrix}}, is outlined in Algorithm~2.
%%%%%%%%%%%%%%%%%%%%%%%%%%%%%%%%%%%%%%%%%%%%%%%%%%%
\begin{algorithm}\caption{: Procedure {\textsc{UpdateMatrix}}$(\psi_\ell(A/2^s), \varphi_\ell(A/2^{s-1}), \ldots,$$\varphi_\ell(A))$.  Under Assumption \ref{ass1}, given $s, \psi_\ell(A/2^s)$  and $\varphi_\ell(A/2^{s-1}),$ $\ldots,  \varphi_\ell(A)$, 
    this procedure returns the approximation of $\psi_\ell(A)$}
    \label{algorithm3}
    \begin{algorithmic}[1] 
        \STATE Initialization: $X_0\gets \psi_\ell(A/2^s);$
            \FOR {$i=s\colon -1\colon 1$} 
            \STATE $X_0 \gets$  {\textsc{NewtonSchulz}$(\varphi_\ell(A/2^{i-1}), X_0)$;} 
            \ENDFOR
              \STATE \textbf{return} $X_0\approx\psi_\ell(A)$
    \end{algorithmic}
\end{algorithm}
%%%%%%%%%%%%%%%%%%%%%%%%%%%%%%%%%%%%%%%%%%%%%%%%%%%
This procedure requires  $\psi_\ell(A/2^s)$  and $\varphi_\ell(A/2^{s-j}),$ $1\leq j\leq s$,  as input.  The computation of $\varphi_\ell(A/2^{s-j}),$ $1\leq j\leq s$, is carried out using the {\tt phipade} algorithm. When  $\ell=1$,  the evaluation of  $\psi_\ell(A/2^s)$  is performed efficiently by {\tt psi1eval}. For $\ell>1$,  $\psi_\ell(z)$ at $z=A/2^s$ 
can be determined by exchanging the numerator and denominator
 of $\widehat R_{d,\ell}(z)$. However,  this approach can be prone to numerical difficulties, since the overall accuracy depends on the conditioning of the numerator of $\widehat R_{d,\ell}(z)$ rather than the denominator. An alternative and more accurate scheme  can be derived  using the results from \cite{gem2023}.  Suppose that  the scaling parameter $s$  in  \eqref{scriterio} is chosen so that  the eigenvalues of $B=A/2^s$ lie within the strip   $\mathbb R\times \mc [-1, 1]$  in the complex plane. Hence, from  Proposition~2 in \cite{gem2023},  we deduce that the {\textsc{NewtonSchulz}} procedure, applied to $\varphi_2(B)$ and $ \psi_1(B),$ quadratically converges to $\psi_2(B)$. This process can be iterated to provide a recursive procedure for computing $\psi_\ell(B)$.   Recall that $\varphi_1(B), \ldots,  \varphi_{\ell-1}(B)$ are generated by {\tt phipade}  for computing  $\varphi_\ell(B)$.

%%%%%%%%%%%%%%%%%%%%%%%%%%%%%%%%%%%%%%%%%%%%%%%%%%%%%%%%%
\subsection{Convergence analysis} \label{convan}

We now demonstrate that, under Assumption~\ref{ass1}, the Newton-Schulz method employed in Procedure {\textsc{UpdateMatrix}} is convergent.

First, we observe that the convergence condition for the Newton-Schulz iteration, when applied to a matrix function, reduces to a scalar inequality.  In particular,  all eigenvalues of $(I-\psi_{\ell}(A) \varphi_{\ell}(2 A))$ must have magnitudes less than $1.$ This condition is equivalent to the following
\begin{equation} \label{rapp_conv}
    \left| \mathcal{C}_\ell(z) \right| := \left|1- \frac{\varphi_\ell (2z)}{\varphi_\ell(z )} \right| <1, \qquad \forall \, z \in \sigma(A), 
\end{equation}
where $\sigma(A)$ denotes the spectrum of $A$.  According  to Assumption~\ref{ass1}, the spectrum of $A$ lies in the  open left half-plane. Consequently, the condition in \eqref{rapp_conv} is satisfied if $\left| \mathcal{C}_\ell(z) \right|<1$ for all complex $z$ with negative real part. The first result concerns real eigenvalues, for which we can obtain a sharper bound.
\begin{theorem}\label{theorel}
Let $z$ be a non-positive real number. Then, for each $\ell \ge 1$
    \begin{eqnarray} \label{eq:modul}
    %\left|1- \frac{\varphi_\ell (2z)}{\varphi_\ell(z)} \right| 
    \left| \mathcal{C}_\ell(z) \right| 
    <\frac{1}{2}.
\end{eqnarray}
\end{theorem}

\begin{proof}
Let $g_\ell(z)$ be the function defined by
\[
g_\ell(z):=  \frac{\varphi_\ell (2z)}{\varphi_\ell(z)}=\frac{   {}_{1}F_{1}(1;\ell+1;2z)}{  {}_{1}F_{1}(1;\ell+1;z)},
\]
where ${}_{1}F_{1}$ denotes the confluent hypergeometric function, as described in \eqref{conflphi}. From \eqref{integrale},  it immediately follows that 
\[
0\leq g_\ell(z)\leq 1, \quad \forall \, z \leq 0.
\]
From \cite[Theorem 2]{MM} we know that for any negative $z$ there exists $\vartheta\in (1,2)$ such that 
\[
{}_{1}F_{1}(1;\ell+1;z)= \displaystyle\frac{1}{1 -\displaystyle\frac{z}{\ell-1 +\vartheta}}.
\]
Hence, it follows that 
\[
g_\ell(z)=\displaystyle\frac{
1 -\displaystyle\frac{z}{\ell-1 +\vartheta_1}}{1 -\displaystyle\frac{2z}{\ell-1 +\vartheta_2}}\geq 
\displaystyle\frac{
1 -\displaystyle\frac{z}{\ell+1}}{1 -\displaystyle\frac{2z}{\ell}}=\frac{\ell}{\ell+1}\frac{\ell+1-z}{\ell-2z}:=r_\ell(z).
\]
Notice that $r_\ell(0)=1$, $\lim\limits_{z\rightarrow -\infty}r_\ell(z)=\displaystyle\frac{1}{2}$. Moreover,  we have $r_\ell'(z)>0$ and, therefore, 
\[
g_\ell(z)\geq r_\ell(z)>\displaystyle\frac{1}{2}, 
\]
which concludes the proof. 
\end{proof}

The case where $z\in \mathbb C^-$ is more involved.  However, a direct proof can be obtained for the case $\ell=1,$  which is addressed in the following theorem.
\begin{theorem}
Let  $z=x+ \mc y\in \mathbb C$ with $x<0$. Then,
    \begin{eqnarray} \label{eq:modul1}
    %\left|1- \frac{\varphi_1 (2z)}{\varphi_1(z)} \right| 
    \left| \mathcal{C}_1(z) \right| <1.
\end{eqnarray}
\end{theorem}

\begin{figure}
  \centering
\subfloat[$\ell=2$]{\includegraphics[width=0.4\textwidth]{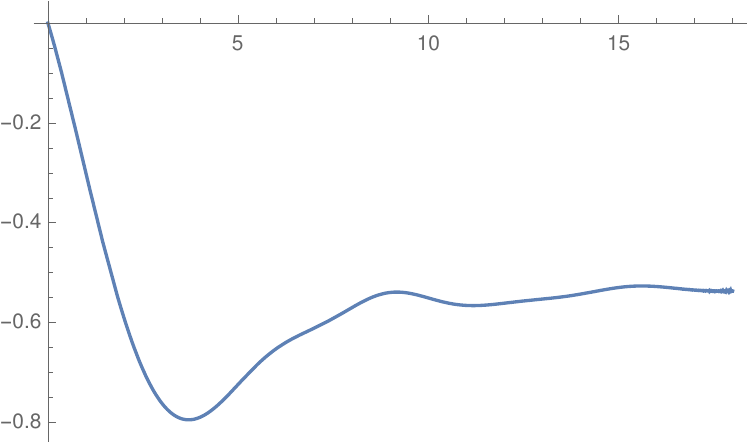}}
\hfill
\subfloat[$\ell=3$]{\includegraphics[width=0.4\textwidth]{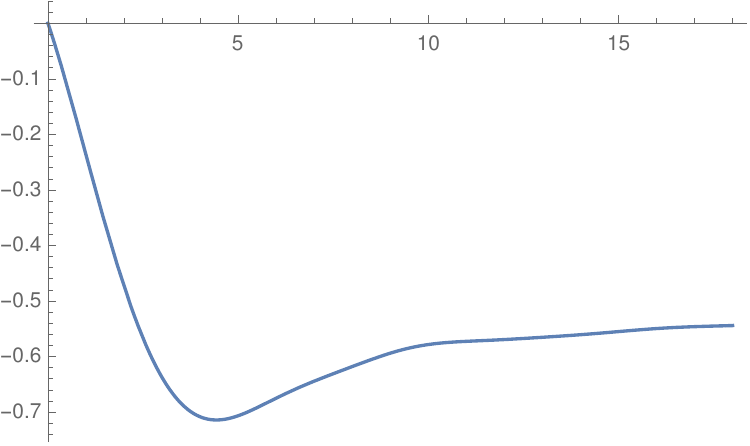}}
\hfill
\subfloat[$\ell=4$]{\includegraphics[width=0.4\textwidth]{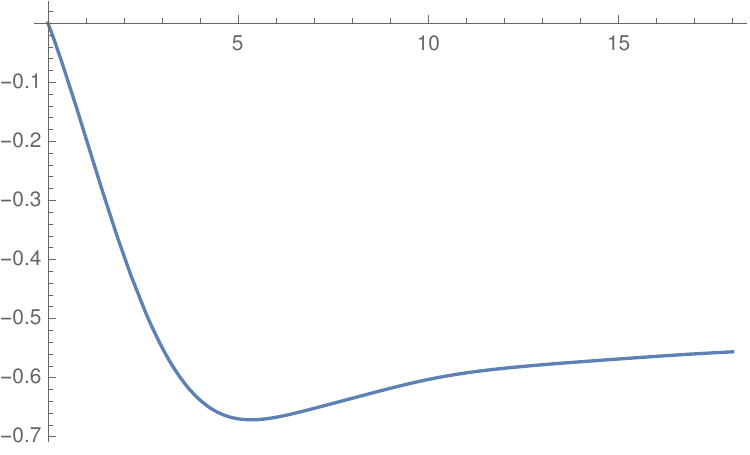}}
\hfill
\subfloat[$\ell=16$]{\includegraphics[width=0.4\textwidth]{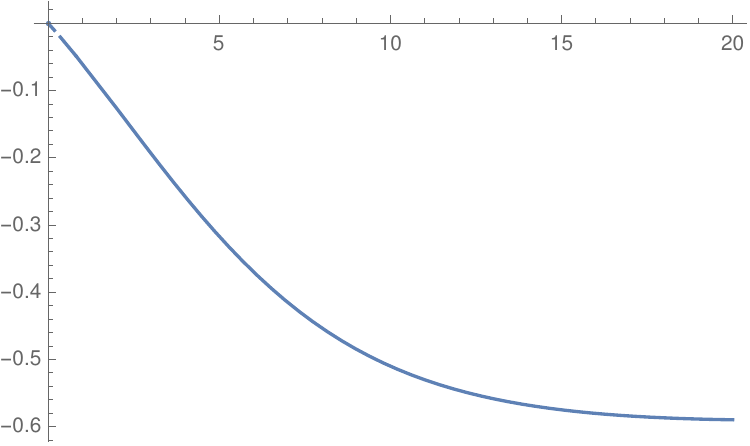}}
  \caption{Plots of $h_\ell(z)$ for $\ell= 2,3, 4, 16,$  with $z\geq 0$.}
\label{f1}
\end{figure}

\begin{proof}
From \eqref{squarfi}  we obtain the relation
\begin{eqnarray*}
 \frac{\varphi_1(2z)}{\varphi_1(z)} &=&   \frac{1}{2} \bigg[ \varphi_0(z)   +   1 \bigg],
\end{eqnarray*}
which implies
\[
\left| \mathcal{C}_1(z) \right| =\left| 1-\frac{\varphi_1(2z)}{\varphi_1(z)} \right| = \frac{1}{2} \left|  1-\varphi_0(z)  \right|.
\]
Therefore, \eqref{eq:modul1} is certainly fulfilled if and only if
\begin{equation} \label{eq:dis}
    \frac{1}{4} \left| 1-\varphi_0(z) \right|^2  <1.
\end{equation}
Now, since  $ z=x+\mc y,$ we can compute 
\begin{eqnarray*}
    1-\varphi_0(z)  &=&   1- e^{x+\mc y} = 1-e^x e^{\mc y} \\
     &=&    1-e^x (\cos y + \mc \sin y) = (1-e^x \cos y) + \mc (-e^x \sin y).
\end{eqnarray*}
Consequently, 
\begin{eqnarray*}
 \left| 1-\varphi_0(z) \right|^2 &=&(1-e^x \cos y)^2 +  (-e^x \sin y)^2 \\
 &=& 1+e^{2x} \cos^2 y-2 e^x \cos y + e^{2x} \sin^2 y \\
 &=& 1+e^{2x} -2 e^x \cos y.
\end{eqnarray*}
Let $w=e^x$, so the inequality in \eqref{eq:dis} becomes 
\[
w^2 -2 \cos y \, w -3<0. 
\]
The roots of the corresponding quadratic equation are
\[
w_1 = \cos y - \sqrt{\cos^2 y +3}, \qquad w_2 = \cos y +\sqrt{\cos^2 y +3} ,
\]
and, therefore, 
the inequality \eqref{eq:dis} holds if and only if
\[
w_1 < e^x < w_2.
\]
Since  $w_1<0$  and $w_2>1,$  we conclude that  this  condition is satisfied when $x<0.$
\end{proof}
For the general case $\ell \geq 2,$  the convergence proof is not fully complete in a theoretical sense; however, we present a partial proof along with numerical evidence to support our findings. First, let us note that
%the  convergence proof  is not  theoretically complete, but we have a partial proof and  numerical evidence. Let us  first note that  
\[
{}_{1}F_{1}(1;\ell+1;-z)=p_\ell(z^2) -z  q_\ell(z^2) 
\]
for suitable power series $p_\ell(z)$ and $q_\ell(z)$.   This  implies that 
\[
{}_{1}F_{1}(1;\ell+1;\mc z)=p_\ell(-z^2)  +\mc z q_\ell(-z^2).
\]
% These two  expressions  establish a  relationship between 
% $\left|1- \displaystyle\frac{\varphi_\ell (-2z)}{\varphi_\ell(-z)} \right|$  and 
% $\left|1- \displaystyle\frac{\varphi_\ell (2\mc z)}{\varphi_\ell(\mc z)} \right|$ for $z\in \mathbb R^+$. 
% Numerically, it is found that 
% \begin{equation}\label{vero}
% h_\ell(z)=\left|1- \frac{\varphi_\ell (2\mc z)}{\varphi_\ell(\mc z)} \right| - 2 \left|1- \frac{\varphi_\ell (-2z)}{\varphi_\ell(-z)} \right|\leq 0, \quad \forall z\geq 0.
% \end{equation}
From \eqref{rapp_conv} these two  expressions  establish a  relationship between $\left| \mathcal{C}_\ell(-z) \right|$  and $\left| \mathcal{C}_\ell(\mc z) \right|$ for $z\in \mathbb R^+$. 
Numerically, it is found that 
\begin{equation}\label{vero}
h_\ell(z):= \left| \mathcal{C}_\ell(\mc z) \right|   - 2\,  \left| \mathcal{C}_\ell(- z) \right|\leq 0, \quad \forall z\geq 0.
\end{equation}
% \begin{equation}\label{vero}
% h_\ell(z)=\left|1- \frac{\varphi_\ell (2\mc z)}{\varphi_\ell(\mc z)} \right| - 2 \left|1- \frac{\varphi_\ell (-2z)}{\varphi_\ell(-z)} \right|\leq 0, \quad \forall z\geq 0.
% \end{equation}
In   Figure  \ref{f1} we show the plot of $h_\ell(z)$ for $\ell=2, 3, 4, 16$.  The limiting behaviour of these curves is in accordance with the asymptotic estimate 
\[
{}_{1}F_{1}(1;\ell+1;z)=\ell (-z)^{-1} (1 + \mathcal{O}(|z|^{-1})),  \quad  z\in \mathbb C^{-}, 
\]
cf. \cite[Eq. 13.1.5]{Abramowitz}, 
which implies 
\begin{equation}\label{limit}
\lim_{z\in \mathbb C^{-}, |z|\rightarrow +\infty} \left| \mathcal{C}_\ell(z) \right|=\displaystyle\frac{1}{2}.
\end{equation}
Assuming that the inequality in \eqref{vero} holds, we can prove the following result.
%then we are able to prove the following. 
\begin{theorem}
 Suppose that \eqref{vero} holds. %is true. 
 Let   $z=x+ \mc y\in \mathbb  C$ with $x<0$. Then,  for any $\ell\geq 2,$ we have $ \left| \mathcal{C}_\ell(z) \right|<1.$ 
\end{theorem}
\begin{proof}
From \cite[Proposition 3.3]{BMG}, we know that the zeros of ${\varphi_\ell (z)}$ are located in the  open 
right half-plane. This means that $\mathcal{C}_\ell(z)$ is an entire function in the open left half-plane  and continuous on the imaginary axis.   Moreover, since $|\mathcal{C}_\ell(z)|=|\mathcal{C}_\ell(\bar z)|,$  we can restrict  the analysis to a quarter of the plane. Let us consider a large square $\Omega =[-x, 0]\times \mc [0, x]$, where $x\gg 0$. By the Maximum Modulus Theorem \cite[Theorem 12-12']{ALV}, it follows that the maximum $\mathtt{m}$ of $|\mathcal{C}_\ell(z)|$ is attained on the closure of $\Omega,$ and we have $|\mathcal{C}_\ell(z)|<\mathtt{m}$ for all $z$ in the  interior of  $\Omega$.  From  Theorem \ref{theorel},  we  have that $|\mathcal{C}_\ell(z)|< 1/2 $ for $z\in [-x, 0]$. 
By using \eqref{limit} we get that $|\mathcal{C}_\ell(z)|<1$ for $z=-x+\mc y$ and $z=y + \mc x$, $0\leq y\leq x$.
Finally,  from  \eqref{vero} and Theorem \ref{theorel}  we obtain that 
\[
\left|\mathcal{C}_\ell(\mc z) \right| \leq 2 \left|\mathcal{C}_\ell(-z) \right| < 1, \quad \forall \, z\in [0, x].
\]
Therefore, we conclude that $\mathtt{m}<1$.
\end{proof}

In the next section, we present numerical applications to assess the robustness and efficiency of the proposed approach.
%%%%%%%%%%%%%%%%%%%%%%%%%%%%%%%%%%%%%%%%%%%%%%%%%%%%
%%%%%%%%%%%%%%%%%%%%%%%%%%%%%%%%%%%%%%%%%%%%%%%%%%%%
\section{Numerical Results}\label{numres}

We performed numerical experiments using matrices from diverse sources, with a particular emphasis on computing $\psi_2(A)$  using Algorithm~2. We focused on several matrices $A\in \mathbb R^{n\times n}$ that arise from parabolic differential problems with non-local boundary conditions. Our test suite consists of the following two matrices:
\begin{enumerate}
 \item[$\bullet$]  
 The first matrix, $\mathcal{A}_1,$ is a diagonal scaling of the standard second derivative matrix  $T=h^{-2} \tridiag[1, -2, 1] \in \mathbb{R}^{n\times n}$, where   $h= 2/(n+1)$.  The matrix is given by
 \[
\mathcal{A}_1= D^{-1}T= \frac{1}{h^2}
\left[\begin{array}{cccc}
\cos(x_1) \\ 
& \cos(x_2)  \\  & & \ddots  \\ & & & \cos(x_n)
\end{array}\right]^{-1}
\left[\begin{array}{cccc}-2 & 1\\ 1& \ddots & \ddots \\& \ddots & \ddots & 1\\ & & 1& -2\end{array}\right].
  \]
 This matrix arises from a discretization on equispaced nodes of the classical heat equation 
 \[
 \displaystyle\frac{\partial u(x,t)}{\partial t}=\frac{1}{\cos(x)}\displaystyle\frac{\partial^2 u(x,t)}{\partial x^2} , \qquad -1\leq x\leq 1, \, 0\leq t\leq 1,
 \]
cf. \cite{Su}.   Observe that $\mathcal{A}_1$ is similar to 
$D^{1/2} \mathcal{A}_1 D^{-1/2}$  which is congruent to the second derivative matrix  $T$. Hence, $\mathcal{A}_1$ has real eigenvalues with negative real part. \label{itemone}

\item[$\bullet$] 
The second matrix,  $\mathcal{A}_2,$ originates from the spatial finite difference discretization of the following advection-diffusion model \cite{LZZ}
\[
\displaystyle\frac{\partial u(x,y,t)}{\partial t}= \displaystyle\frac{\partial^2 u(x,y,t)}{\partial x^2} + \displaystyle\frac{\partial^2 u(x,y,t)}{\partial y^2} -10x \displaystyle\frac{\partial u(x,y,t)}{\partial x} -100 y \displaystyle\frac{\partial u(x,y,t)}{\partial y}
\]
with $0\leq x,y \leq 1$, $0\leq t\leq 1$. The matrix $\mathcal{A}_2$ is  generated directly by  means of the MATLAB function {\tt fdm}$\_${\tt 2d}$\_${\tt matrix} from LYAPACK toolbox \cite{VM}. \label{itemtwo}
\end{enumerate}
To benchmark the performance of  Algorithm~2, we have conducted numerical experiments using MATLAB (version R2019a).  
The algorithm {\tt{phipade}} from the EXPINT package \cite{Skaflestad,EXPINT} is used to compute the input data.
The value of $s$  is  determined  by {\tt{phipade}} according to \eqref{scriterio} with $\theta=4$. 
The squaring process for the matrix $\varphi$-functions is  also performed by {\tt{phipade}}.  The seeds of this process, i.e., $\varphi_j(A/2^s)$ for $0\leq j\leq 2$, are computed using the corresponding Padé approximants with  $d=7$.  To generate the matrix $\psi_2(A/2^s)$, we apply Procedure {\textsc{NewtonSchulz}$(\varphi_2(A/2^s), \psi_1(A/2^s))$}, where $\psi_1(A/2^s)$ is returned by {\tt psi1eval} with $d=7$. \\

\begin{figure}
  \centering
{\includegraphics[width=0.6\textwidth]{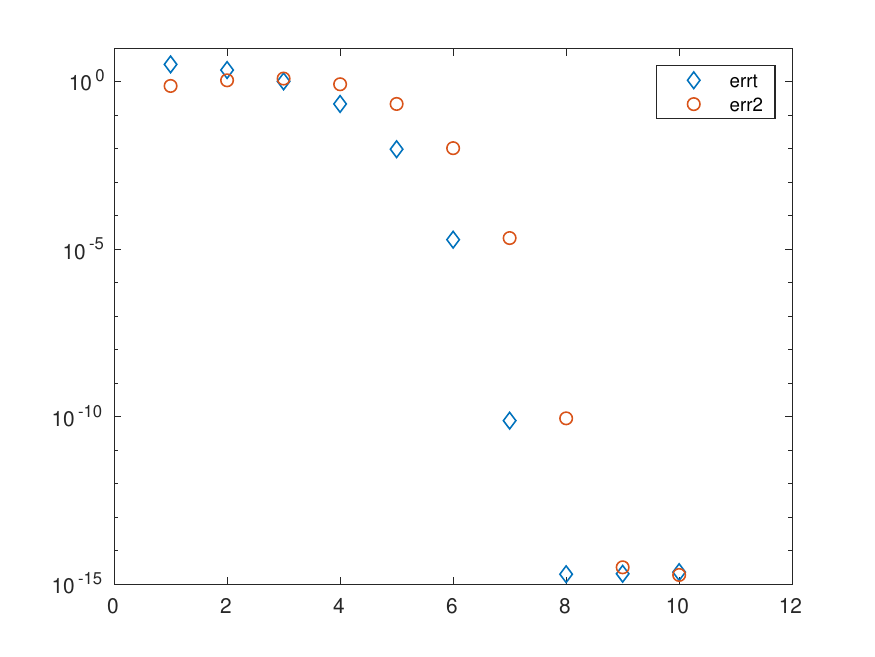}}
  \caption{Convergence of  the Newton-Schulz  algorithm  applied to $\varphi_2(\mathcal{A}_1/2^{19})$ with starting point $\psi_1(\mathcal{A}_1/2^{19})$. The plots illustrate the convergence of $err_2$ of  Algorithm~1  for the matrix $\mathcal{A}_1$ with $n=1024$. }
  \label{updateorderfig}
\end{figure}

\begin{figure}
  \centering
{\includegraphics[width=0.6\textwidth]{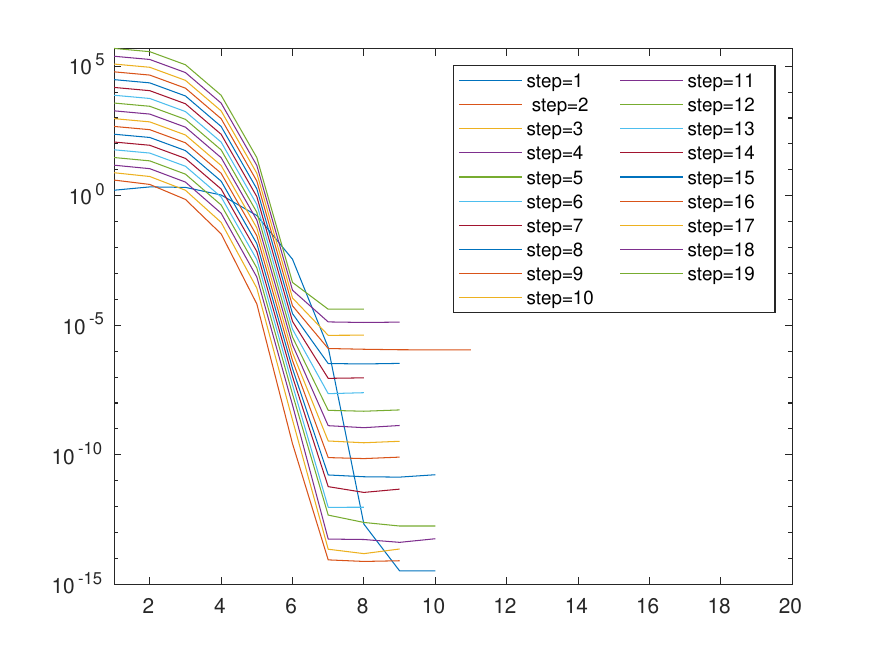}}
  \caption{Convergence of  the Newton-Schulz  algorithm  called  by the {\textsc{UpdateMatrix}} procedure to perform the squaring steps for matrix $\psi$-functions. The plots illustrate the convergence of $err_2$ of  Algorithm~1  in each squaring step for the matrix $\mathcal{A}_1$ with $n=1024$. }
  \label{conv_history}
\end{figure}

For  the matrix $\mathcal{A}_1$  with  $n=1024$ {\tt{phipade}}  finds $s=19$. In Table \ref{t1} we show the absolute error   $err=\| \psi_1(A/2^s)- \mathcal{R}_{d,1}(A/2^s)\|_\infty$, where $\mathcal{R}_{d,1}(A)$ is the approximation provided by {\tt psi1eval}
 and $\psi_1(A/2^s)$  is  computed using the built-in function  {\tt expm} and the backslash operator. 
 Additionally, we provide a numerical estimate 
$err_{est}$ of the error upper bound derived from Theorem~\ref{phipadetheo}. This estimate is obtained by combining the upper bound from \eqref{bb} with a lower bound on $|\mathcal D_{d,1}(z)|$ for $|z|\leq 4.$
\begin{table}
\begin{tabular}{|l|l|l|l|l|}  \hline
			$d$ & $6$ &$ 7$ & $8$ & $9$ \\ 	\hline
			$err$ & $7.9 \cdot 10^{-8}$ & $1.1 \cdot  10^{-9}$ & $9.7\cdot  10^{-11}$& $9.7\cdot  10^{-11}$ \\ \hline
            $err_{est}$ & $1.9 \cdot  10^{-6}$& $3.4 \cdot  10^{-8}$& $4.8 \cdot  10^{-10}$& $5.4 \cdot  10^{-12}$ \\ \hline
			\end{tabular}
       \caption{Comparison of computed errors and theoretical estimates for the $[d/d]$-Pad\'e approximant of $\psi_1(\mathcal{A}_1/2^{19})$. \label{t1}} 
\end{table}
In Figure~\ref{updateorderfig}, we illustrate the behavior of the Newton-Schulz iteration for computing  $\psi_2(A/2^s)$ starting from the approximation of $\psi_1(A/2^s)$ returned by {\tt psi1eval}. For comparison, we plot the error term 
 $err_2$  in  Algorithm~1 versus  the reference error $err_t=\parallel X_k-(\varphi_2(A/2^s))^{-1} \parallel_\infty$ where the inverse is computed using the backslash operator.
In Figure~\ref{conv_history} we show the convergence history of the Newton-Schulz  iterations performed by the {\textsc{UpdateMatrix}} procedure.
\\

\begin{figure}
  \centering
{\includegraphics[width=0.6\textwidth]{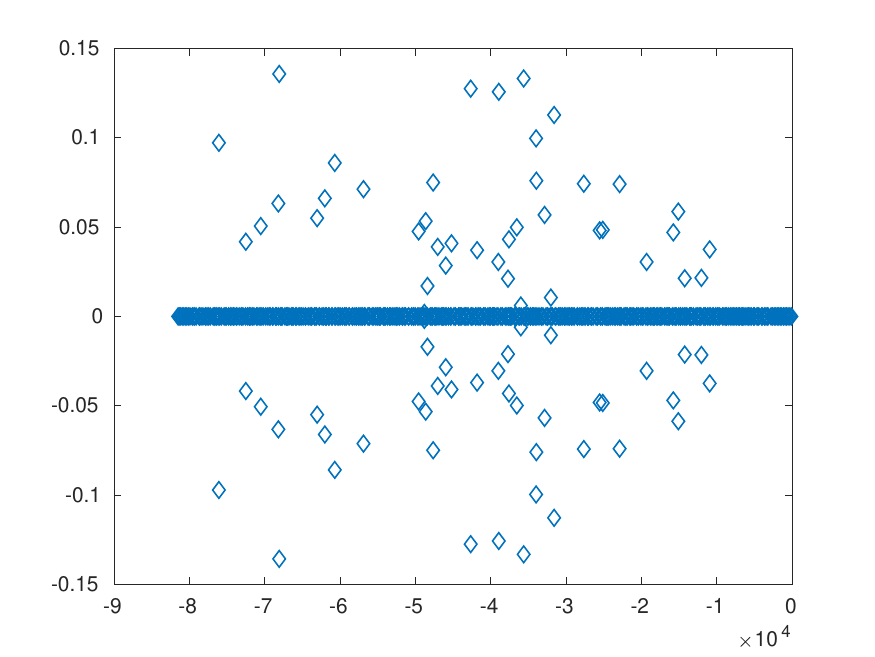}}
  \caption{Illustration of the spectrum of $\mathcal{A}_2$ for $n=10000$.}
  \label{eigenv}
\end{figure}

\begin{figure}
  \centering
{\includegraphics[width=0.6\textwidth]{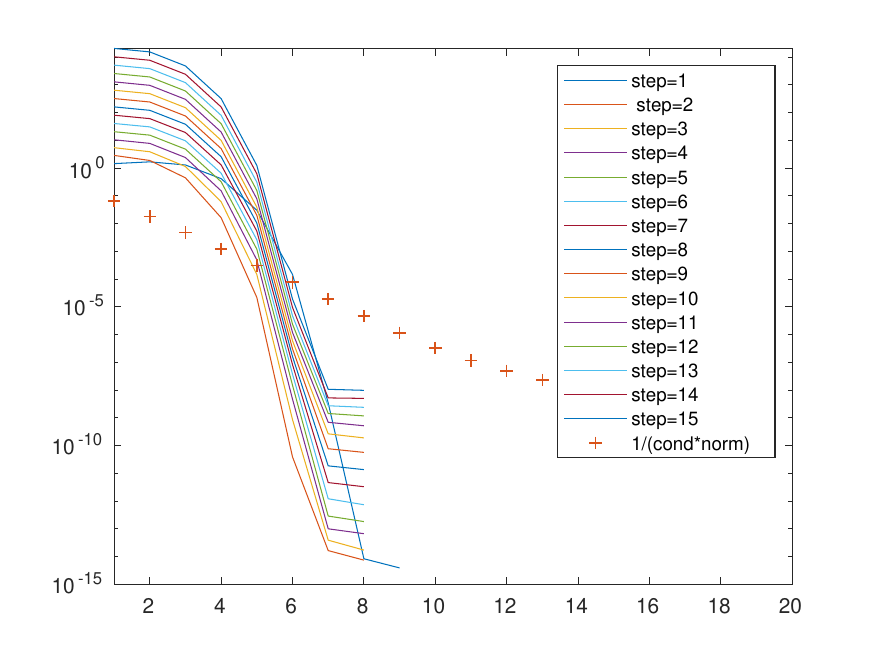}}
  \caption{Convergence of  the Newton-Schulz  algorithm  called  by the {\textsc{UpdateMatrix}} procedure to perform the squaring steps for matrix $\psi$-functions. The plots illustrate the convergence of $err_2$ of  Algorithm 1  in each squaring step for the matrix $\mathcal{A}_2$  with $n=10000$. }
  \label{conv_history2}
\end{figure}

The matrix $\mathcal{A}_2$ is unsymmetric.  For $n=10000$ 
its eigenvalues are found to lie in the open left half-plane. In Figure~\ref{eigenv}, we show the spectrum of 
$\mathcal{A}_2,$ computed using the built-in function {\tt eig}.
In Figure~\ref{conv_history2} we show the convergence history of the Newton-Schulz  iterations performed by the {\textsc{UpdateMatrix}} procedure applied to the matrix $\mathcal{A}_2$  with $n=10000$. For comparison, we also show the plot of the reciprocal of the product of the condition number and the norm of the matrix, evaluated at each step. This clearly highlights that the iteration halts when the approximation reaches the limits defined by the conditioning estimates.

%%%%%%%%%%%%%%%%%%%%%%%%%%%%%%%%%%%%%%%%%%%%%%%%%%%%%
\section{Conclusions}\label{sec:fine}

In this paper, we have introduced a scaling-and-squaring method for evaluating the  matrix $\psi$-functions, which are  the inverses of matrix $\varphi$-functions.  The method is based on the Newton-Schulz  algorithm  for matrix inversion, which is  advantageous over other methods for high-performance computing because it is rich in matrix-matrix multiplications. 
Furthermore, this approach can be related to Krylov-type methods such as GMRES, allowing for adjustments that make it suitable for computing the action of matrix $\psi$-functions on a vector. 
Future work will focus on extending this method in that direction. Additionally, we plan to implement our approach in a parallel distributed computing environment to further enhance its scalability and efficiency.

%%%%%%%%%%%%%%%%%%%%%%%%%%%%%%%%%%%%%%%%%%
%%%%%%%%%%%%%%%%%%%%%%%%%%%%%%%%%%%%%%%%%
 \section*{Acknowledgements}
The authors are members of the INdAM research group GNCS.

%%%%%%%%%%%%%%%%%%%%%%%%%%%%%%%%%%%%%%%%%%%%%%%%%%%%%
\section*{Funding}
Lidia Aceto is partially supported by the “INdAM - GNCS Project", with code CUP$\_$E53C23001670001. Luca Gemignani is  partially supported by  European Union - NextGenerationEU under the National Recovery and Resilience Plan (PNRR) - Mission 4 Education and research - Component 2 From research to business - Investment 1.1 Notice Prin 2022 - DD N. 104  2/2/2022, titled Low-rank Structures and Numerical Methods in Matrix and Tensor Computations and their Application, proposal code 20227PCCKZ – CUP I53D23002280006  and by the Spoke 1 “FutureHPC \& BigData”  of the Italian Research Center on High-Performance Computing, Big Data and Quantum Computing (ICSC)  funded by MUR Missione 4 Componente 2 Investimento 1.4: Potenziamento strutture di ricerca e creazione di "campioni nazionali di R\&S (M4C2-19)" - Next Generation EU (NGEU).

%%%%%%%%%%%%%%%%%%%%%%%%%%%%%%%%%%%%%%%%%%%%%%%%%%%%%%%%%%%%
\section*{Declarations}
{\bf{Conflict of interest}} The authors declare that they have no conflict of interest.\\
{\bf{Data Availability}} Enquiries about data availability should be addressed to the authors.

%%%%%%%%%%%%%%%%%%%%%%%%%%%%%%%%%%%%%%%%%%%%%%%%%%%%%%%%%%%%
%\bibliographystyle{tfnlm}
%\bibliography{bibliographyPADE}
%\bibliographystyle{sn-mathphys-ay}
%\bibliography{sn-bibliography}
%% BioMed_Central_Bib_Style_v1.01

\end{document}